\documentclass[preprint,12pt]{elsarticle}
\usepackage{amsmath}
\usepackage{amsthm}
\usepackage{amsfonts}
\usepackage{units}
\usepackage{amssymb}
\usepackage{extarrows}
\usepackage{tikz}
\usepackage{verbatim}
\usetikzlibrary{arrows, matrix}
\newcommand{\Z}{\mathbb{Z}}

\renewcommand{\phi}{\varphi}
\newcommand{\im}{\operatorname{im}}

\newcommand{\spec}{\operatorname{Spec}}

\renewcommand{\H}{\mathcal{H}}
\newcommand{\Inv}{\operatorname{Inv}}

\journal{Journal of Pure and Applied Algebra}

\begin{document}

\begin{frontmatter}

\title{Vanishing of Degree 3 Cohomological Invariants}


\author{Rebecca Black\corref{cor1}}
\cortext[cor1]{University of Maryland, College Park}
\ead{rblack1@math.umd.edu}
\address{4176 Campus Drive - William E. Kirwan Hall
College Park, MD 20742-4015}

\begin{abstract}
For a complex algebraic variety $X$, we show that triviality of the sheaf cohomology group $H^0(X,\mathcal{H}^3)$ occurring on the second page of the Bloch-Ogus spectral sequence \cite{bloch-ogus} follows from a condition on the integral Chow group $CH^2X$ and the integral cohomology group $H^3(X, Z)$. In the case that $X$ is an appropriate approximation to the classifying stack $BG$ of a finite $p$-group $G$, this result states that the group $G$ has no degree three cohomological invariants. As a corollary we show that the nonabelian groups of order $p^3$ for odd prime $p$ have no degree three cohomological invariants.

\end{abstract}

\begin{keyword}
algebraic geometry \sep cohomological invariant \sep motivic cohomology \sep Chow ring

\MSC[2010] 14C15  \sep 14F42  \sep 20G10  

\end{keyword}

\end{frontmatter}

\section{Introduction} Let $G$ be a finite $p$-group of order $p^n$, considered as an algebraic group over $\mathbb{C}$.  In this paper we employ the tools of the Bloch-Ogus spectral sequence and the motivic cohomology ring of the classifying space $BG$ in order to examine in detail the relationship between the Chow ring of $G$ and the ring of cohomological invariants of $G$ in low degree.  In particular, our main result is that if the cycle class map 
$$cl: CH^2G \to H^4(BG, \Z)$$ is an isomorphism, then there are no non-trivial degree three cohomological invariants of $G$.  There has been a lot of progress recently in computing the Chow rings of various classes of $p$-groups, so we know that we have this isomorphism in certain cases. (See for example \cite{totarobook} for an excellent overview of recent progress.) 

Ideally similar techniques could be employed to explicitly relate the Chow ring to vanishing of invariants in higher degree as well, but a more detailed computational understanding of the motivic cohomology of $G$ is necessary to extend this method to higher degrees.

\section{Chow groups and cohomological invariants}

Totaro defines the Chow groups of the group $G$ in terms of finite-dimensional approximations to the classifying stack $BG$.  Suppose that $V$ is a representation of $G$, and let $S \subset V$ be the locus on which the stabilizers are non-trivial, with $\operatorname{codim}S = d$.  Let $X = (V-S)/G$ be the quotient variety.  Then 
$$CH^iG = CH^iX \,\, \text{for } i < d \text{.}$$
Totaro proved the existence of a sequence of such representations $V_n$ with the codimension of $S_n$ going to infinity; see \cite{totarochow} for a good exposition. Throughout this paper we will freely assume that we have a variety $X$ of this form where we have taken the representation to be of high enough dimension that $X$ has the same invariants and cohomology as $G$ in low dimensions.

A cohomological invariant of $G$ is a natural transformation of functors 
$$\eta: H^1(-, G) \to H^*(-, \Z/p) \text{,}$$
where $H^1(K,G)$ is the first nonabelian Galois cohomology set (which can be thought of as isomorphism classes of $G$-torsors over $K$), and $H^*(K,\Z/p)$ is the abelian Galois cohomology ring.  For our purposes, however, this is not the most convenient way to think of cohomological invariants.  Given a quotient variety $X = (V-S)/G$ as above, with $\operatorname{codim}S \geq 2$, the generic fiber $T$ of the map $V-S \to X$ is a versal $G$-torsor, meaning that any given cohomological invariant is actually completely defined by its value on that specific torsor (see discussion in \cite{guillot}). Since $T$ is defined over $\spec k(X)$, its image under an invariant $\eta$ will lie in the Galois cohomology group $H^d(k(X), Z/p)$ for some degree $d$.  Hence we can identify the group of degree $d$ cohomological invariants of $G$ with a certain subset of $H^d(k(X), \Z/p)$.

In fact, we can say much more about that certain subset: Given a point $x \in X$ with $\operatorname{codim}\overline{\{x\}} = 1$, we get a residue map $$\nu_x: H^d(k(X),\Z/p) \to H^{d-1}(k(x), \Z/p) \text{,}$$
where $k(x)$ is the stalk at $x$. If a class $\eta_T \in H^d(k(X),\Z/p)$ is the image of a versal torsor under an invariant, then $\nu_x(\eta_T) = 0$ for all such $x$; conversely, Totaro shows that if $\operatorname{codim}S \geq 2$, every class in the kernel of $\nu_x$ for all $x$ does in fact define a cohomological invariant (letter to Serre, reprinted in \cite{GMS}). Therefore we have the identification
$$\Inv^dG = \ker\left( H^d(k(X),\Z/p) \to \coprod_{x \in X^{(1)}} H^{d-1}(k(x),\Z/p) \right) \text{,}$$
where $x \in X^{(1)}$ ranges over all codimension one points.

\section{Bloch-Ogus spectral sequence and stable cohomology}
In their 1974 paper, Bloch and Ogus showed that the product of residue maps considered above is part of a flasque resolution of the sheaf $\H^d$ on $X$, defined as the sheafification of the Zariski presheaf $U \mapsto H^d_{\text{\'et}}(U, \Z/p)$.  Therefore we can actually think of the kernel as a sheaf cohomology group, and we get 
$$\Inv^dG = H^0(X, \H^d) \text{.}$$
This sheaf cohomology group appears as the $E_2^{0,d}$ term of the Bloch-Ogus spectral sequence for $X$, which converges to the \'etale cohomology $H^*_{\text{\'et}}(X, \Z/p)$.  With our assumptions on $X$ and the base field $k$, we can in fact identify these \'etale cohomology groups with the group cohomology $H^*(G, \Z/p)$ in low degree.

The diagonal entries $E_2^{r,r}$ are isomorphic to the mod $p$ Chow groups $CH^rX \otimes \Z/p \cong CH^rG \otimes \Z/p$.  Hence the differential $\delta: E_2^{0,3} \to E_2^{2,2}$ combined with the maps to and from the abutment give an exact sequence:
$$H^3(G,\Z/p) \to \Inv^3G \xrightarrow{\delta} CH^2G \otimes \Z/p \to H^4(G,\Z/p) \text{.}$$

Our basic plan of attack is to show that both the kernel and the image of $\delta$ are trivial, which forces $\Inv^3G = 0$.  Triviality of the image will follow immediately from the assumption on the integral Chow groups, since the map to $H^4(G,\Z/p)$ coincides with the mod $p$ cycle class map.  We will have to work a bit harder to show that the kernel is trivial.  The kernel of $\delta$ is precisely the classes that survive to the group cohomology, also known as the stable cohomology (as discussed in \cite{bog-boh}).  Hence the vanishing of this kernel is equivalent to the vanishing of degree three stable cohomology, which is shown for several cases of $p$-groups in \cite{bog-boh}.  We will use a somewhat different argument that makes use of the relationship of the sheaf cohomology groups $H^r(X,\H^s)$ with the motivic cohomology ring $H^{*,*'}(X,\Z/p)$.

\section{Motivic cohomology}

This section summarizes a few important properties of the motivic cohomology ring $H^{*,*'}(X,\Z/p)$ associated to a variety $X$. We do not attempt a complete discussion of the definition of this ring here; see for example \cite{voevodsky} for details.  The beauty of the motivic cohomology ring for us is that it specializes for certain indices to both the mod $p$ Chow groups and the \'etale cohomology groups.  Specifically, Voevodsky and others have shown the following, for a projective variety $X$:

$$H^{m,n}(X,\Z/p) \cong \left\{ \begin{array}{ll} 0 & \,\, \text{if $m > 2n$;} \\ CH^nX \otimes \Z/p & \,\, \text{if $m = 2n$;} \\ H^m_{\text{\'et}}(X,\Z/p) & \,\, \text{if $m \leq n$.} \end{array} \right.$$

Let $\tau$ denote a generator of $H^{0,1}(\spec(k),\Z/p) \cong \Z/p$. Then the cup product gives a map $\times \tau: H^{m,n}(X,\Z/p) \to H^{m,n+1}(X,\Z/p)$. Our argument makes use of the following long exact sequence, which relates this map to the sheaf cohomology groups that appear in the Bloch-Ogus spectral sequence (see \cite{yagita}):
\begin{align*} \cdots \rightarrow &H^{m,n-1}(X,\Z/p) \xrightarrow{\times \tau} H^{m,n}(X,\Z/p) \rightarrow \\&H^{m-n}(X, \H^n) \rightarrow H^{m+1,n-1}(X,\Z/p) \xrightarrow{\times \tau} \cdots \end{align*}

Finally, we will also use the fact that there is motivic cohomology with integer coefficients as well, and in particular there is an isomorphism 
$$H^{2n,n}(X,\Z) \cong CH^nX \text{.}$$

\section{Main Theorem}

We are now ready to state and prove our main result on the vanishing of degree three cohomological invariants.

\theoremstyle{plain}
\newtheorem{deg3}{Theorem}[section]
\begin{deg3} Let $X$ be a variety over $\spec \mathbb{C}$ satisfying the following two properties: \begin{itemize}
    \item [(i)] $CH^2X \cong H^4(X,\Z)$;
    \item [(ii)] There is some power $p^n$ with $p^n H^3(X,\Z) = 0$.
\end{itemize}
Then $H^0(X,\mathcal{H}^3) = 0$.  In particular, if $X$ is an approximation of the classifying stack $BG$ for an algebraic group $G$ such that the above two conditions hold, then $\operatorname{Inv}^3G = 0$.
\end{deg3}

\begin{proof} The group $H^0(X,\mathcal{H}^3)$ fits into the following long exact sequence:
\begin{align*} \cdots \to H^{3,2}(X, \Z/p) \xrightarrow{\times \tau} &H^{3,3}(X, \Z/p) \to H^0(X, \mathcal{H}^2) \to \\&H^{4,2}(X, \Z/p) \xrightarrow{\times \tau} H^{4,3}(X, \Z/p) \to \cdots \text{.} \end{align*}
Therefore, we get our result if we can show that \begin{itemize}
    \item [(a)] $\times \tau: H^{4,2}(X, \Z/p) \to H^{4,3}(X, \Z/p)$ is injective.
    \item [(b)] $\times \tau: H^{3,2}(X, \Z/p) \to H^{3,3}(X, \Z/p)$ is surjective, and
\end{itemize}  The injectivity is easier to show, so we will do that first. We know that $H^{4,2}(X, \Z/p) \cong CH^2X \otimes \Z/p$ is the mod $p$ Chow group. The mod $p$ cycle class map $$c: CH^2X \otimes \Z/p \to H^4(X, \Z/p)$$ agrees with the change of coefficients map $H^4(X, \Z) \to H^4(X,\Z/p)$ induced by the short exact sequence $$0 \to \Z \to \Z \to \Z/p \to 0 \text{,}$$ meaning its kernel is exactly $pCH^2X$.  This shows that $c$ is injective.  Since we can identify $c$ with the map $\times \tau^2$ on motivic cohomology, we have shown (a).

For (b), denote by $\beta$ the connecting homomorphism $\beta: H^{3}(X, \Z/p) \to H^4(X, \Z)$; the plan of attack is first to show that $$\operatorname{ker}(\beta) \subseteq \operatorname{im}(\tau) \subseteq H^{3,3}(X,\Z/p) \cong H^3(X,\Z/p) \text{,}$$ and then to show that any class in $H^{3}(X,\Z/p)$ is equivalent to $\operatorname{ker}(\beta)$ mod the image of $\tau$.

The key to the first step is that, for any exponent $n$ the short exact sequence $$0 \to \Z/p \to \Z/p^{n+1} \to \Z/p^n \to 0$$ induces connecting maps on both \'etale cohomology and motivic cohomology:
$$\beta_{\text{mot}}: H^{*,*}(X,\Z/p^n) \to H^{*+1,*}(X,\Z/p) \text{;}$$
$$\beta_{\text{et}}: H^{*}(X,\Z/p^n) \to H^{*+1}(X,\Z/p) \text{.}$$

Under the isomorphisms from the Beilinson-Lichntenbaum conjecture, then, $\beta_{\text{et}}$ maps from $H^{*,*}(X,\Z/p^n)$ to $H^{*+1,*+1}(X, \Z/p)$, and we have that $$\beta_{\text{et}} = \tau \circ \beta_{\text{mot}} \text{.}$$ Therefore, crucially for us, $\operatorname{im}(\beta_{\text{et}}) \subseteq \operatorname{im}(\tau)$. 

Now let $x \in \operatorname{ker}(\beta) \subseteq H^3(X,\Z/p)$. Then we can pull $x$ back to a class $\Tilde{x} \in H^3(X,\Z)$. By assumption, $H^3(X,\Z)$ is $p^n$-torsion for some $n$, meaning that $\Tilde{x}$ in turn comes from a class $\overline{x} \in H^2(X,\Z/p^n)$. Then we have $x = \beta_{\text{et}}(\overline{x}) \in \operatorname{im}(\tau)$ as desired.

For the general case, we now assume that $\beta(x) \neq 0 \in H^4(X,\Z)$. Recall that by assumption $H^4(X,\Z) \cong CH^2X \cong H^{4,2}(X,\Z)$; we write $y \in H^{4,2}(X,\Z)$ for the image of $\beta(x)$ under this isomorphism.  Since $py = 0$, we have $y = \beta(x')$ for some $x' \in H^{3,2}(X,\Z/p)$.  Then $\beta(\tau x') = \beta(x) \in H^4(X,\Z)$ (where we abuse notation a bit by conflating $\tau x'$ and its image under the isomorphism $H^{3,3}(X,\Z/p) \cong H^3(X,\Z/p)$).  By the previous case, then, $\beta(x - \tau x') = 0$, so $x - \tau x' \in \im(\tau)$; therefore we also have $x \in \im(\tau)$ as desired.

\end{proof}

In the case that $X = (V-S)/G$ is an approximation to $BG$ as described above, with $|G| = p^n$, we do automatically have that $H^3(X,\Z) \cong H^3(G,\Z)$ is $p^n$-torsion, so the second condition of the theorem is automatically satisfied. Therefore we have shown that for finite $p$-groups $G$, if the degree two cycle class map is an isomorphism then $G$ has no nontrivial degree three cohomological invariants. For example, Yagita proved that the cycle class map is an isomorphism in all degrees for the two nonabelian groups of order $p^3$ for odd primes $p$ \cite{yagita}, meaning by our result these groups have no cohomological invariants of degree three.

\section{Acknowledgements} Infinite thanks to Patrick Brosnan for his help and patience as I struggled to work through the last details of this proof; and thank you also to Nobuaki Yagita and Masaki Kameko for extremely helpful suggestions and clarifications.

    \bibliographystyle{elsarticle-num}
    \bibliography{cohoinv}

\end{document}